\documentclass[10pt,reqno]{amsart}

\usepackage{dsfont}
\usepackage[bottom]{footmisc}
\usepackage{amssymb, amsmath, amsthm, mathrsfs,amsfonts,latexsym}
\usepackage[utf8]{inputenc}
\usepackage{enumerate,color,comment}
\usepackage{hyperref}
\usepackage{cleveref}

\newtheorem{thm}{Theorem}[section]

\newtheorem{prop}[thm]{Proposition}
\newtheorem{defn}[thm]{Definition}

\numberwithin{equation}{section}
\newcommand {\ds}{\ \mbox{d}s}

\newcommand{\bel}{\begin{equation} \label}
\newcommand{\ee}{\end{equation}}
\def\beq{\begin{equation}}
\def\eeq{\end{equation}}
\newcommand{\bea}{\begin{eqnarray}}
\newcommand{\eea}{\end{eqnarray}}
\newcommand{\beas}{\begin{eqnarray*}}
\newcommand{\eeas}{\end{eqnarray*}}
\newcommand{\pd}{\partial}

\newcommand{\dd}{\mbox{d}}

\newcommand{\ep}{\varepsilon}

\newcommand{\R}{\mathbb{R}}
 
\newcommand{\N}{\mathbb{N}}

\newcommand{\cB}{\mathcal{B}}
\newcommand{\cC}{\mathcal{C}}

\newcommand{\cL}{\mathcal{L}}
\newcommand{\sL}{\mathscr{L}}
\newcommand{\sH}{\mathscr{H}}
\newcommand{\sV}{\mathscr{V}}

\renewcommand{\div}{\mathrm{div}\,}  

\allowdisplaybreaks

\def\epsilon{\varepsilon}
\def\phi {\varphi}

\providecommand{\abs}[1]{\left\lvert#1\right\rvert}
\providecommand{\norm}[1]{\left\lVert#1\right\rVert}

\providecommand{\abs}[1]{\left\lvert#1\right\rvert}
\providecommand{\norm}[1]{\left\lVert#1\right\rVert}

\title[Time-fractional partial differential equations of piecewise constant time-varying order]
{\bf On time-fractional partial differential equations of time-dependent piecewise constant order}

\author{Yavar Kian$^\triangle$}
\address[$\triangle$]{
Univ Rouen Normandie, CNRS, Normandie Univ, LMRS UMR 6085, F-76000 Rouen, France}
\email{yavar.kian@univ-rouen.fr}

\author{Mari\'an Slodi\v{c}ka$^\bigstar$}
\address[$ \bigstar$ ]{Department of Electronics and Information Systems, research group of Numerical Analysis and Mathematical Modeling (NaM$^2$), Ghent University, Krijgslaan 281, S8, Gent 9000, Belgium}
\email{marian.slodicka@ugent.be}

\author{\'Eric Soccorsi$^\Box$}
\address[$\Box$]{Aix-Marseille Univ, Université de Toulon, CNRS, CPT, Marseille, France}
\email{eric.soccorsi@univ-amu.fr}

\author{Karel Van Bockstal$^\spadesuit$}
\address[$\spadesuit$]{Ghent Analysis \& PDE center, Department of Mathematics: Analysis, Logic and Discrete Mathematics, Ghent University, Krijgslaan 281, 9000 Ghent, Belgium}
\email{karel.vanbockstal@ugent.be}

\date{\today}

\begin{document}

\begin{abstract} 
This contribution considers the time-fractional subdiffusion with a time-dependent variable-order fractional operator of order $\beta(t)$. It is assumed that $\beta(t)$ is a piecewise constant function with a finite number of jumps. A proof technique based on the Fourier method and results from constant-order fractional subdiffusion equations has been designed.  
This novel approach results in the well-posedness of the problem. 
\end{abstract}

\maketitle


\section{Introduction}

The idea of fractional calculus (FC) dates back to the 17th century when mathematicians like Leibniz and L'H\^{o}pital pondered the meaning of differentiation and integration of noninteger orders. Significant progress was made in the 19th century thanks to the work of mathematicians like Liouville and Riemann.

Fractional calculus has applications in various scientific and engineering fields, including physics, engineering, signal processing, finance, and more. It has proven to be a powerful tool for describing systems with long-range memory, fractal phenomena, and non-local behaviour. The fractional calculus framework provides a deeper understanding of complex phenomena that integer-order calculus cannot adequately describe.

Machado et al. \cite{MACHADO2011} summarised the historical perspective on the major developments in fractional calculus since the 1970s. 
 Please note that most papers studied constant-order (CO) fractional operators.
 There exist various types of CO fractional derivatives (Caputo, Riemann-Liouville,  Gr\"{u}nwald-Letnikov, etc.), see, e.g. \cite{Podlubny1999}. The Caputo and the Riemann-Liouville variations can be expressed in terms of the Riemann-Liouville kernel $g_{1-\beta}$:  
\[
	g_{1-\beta}(t) = \frac{t^{-\beta}}{\Gamma(1-\beta)},\qquad t>0,\ 0<\beta<1,
\]
where $\Gamma$ represents the Gamma function. 
This kernel obeys
\[
    (-1)^j g^{(j)}_{1-\beta}(t) \ge 0, \quad \forall t\ge 0, j=0,1,2; \quad g^\prime_{1-\beta} \not \equiv 0, 
\]
and, therefore, it is strongly positive definite by \cite[Corollary 2.2]{Nohel1976}. This creates a powerful position when establishing energy estimates for a solution.

In the literature, several definitions of fractional derivatives and integrals in variable-order (VO) FC can be found as generalisations of their CO counterparts, see e.g. \cite{almeida2019,Coimbra2003,Ding2020,FF2012,K2023,KS2023,KSY2018,Lorenzo2002,Sun2009,VanBockstal2021,Zhuang2009,SunChangZhangChen2019}.
One of the possibilities uses the following time-dependent Riemann-Liouville kernel 
\[
    \left( g_{1-\beta(t)}\right)(t) = \frac{t^{-\beta(t)}}{\Gamma(1-\beta(t))}, \qquad t>0,\ 0<\beta(t)<1.
\]
 Then, the Caputo VO fractional derivative reads as follows
 \[
    \left(\partial_t^{\beta(t)} u\right)(t) = \left(g_{1-\beta(t)} \ast \partial_t u\right)(t) = \int_0^t \frac {(t-s)^{-\beta(t-s)}}{\Gamma(1-\beta(t-s))} \partial_t u(s)\ds, \qquad t>0,\ 0<\beta(t)<1.
\]
Another possibility for a slightly different VO derivative is \cite{Coimbra2003,Sun2009}
\[
 \left(\partial_t^{\beta(t)} u\right)(t) =  \frac {1}{\Gamma(1-\beta(t))} \int_0^t  (t-s)^{-\beta(t)}  \partial_t u(s)\ds, \qquad t>0,\ 0<\beta(t)<1.
\] 
This definition will be considered in this paper. Polymers, plastics, rubber and oil are all viscoelastic materials that are used in biology, medicine, chemical engineering and other fields.
CO fractional viscoelasticity models do not account for the evolution of the microstructure of viscoelastic materials during deformation, which changes the mechanical properties of deforming materials. 
On the other hand, in VO fractional models with a fractional time-derivative of time-dependent order, the variation of the fractional order is used for describing the change of the mechanical properties of the material under various loading conditions.
For this reason, they can accurately predict the viscoelastic behaviour of soft materials, see e.g., \cite{Gao2022,Ramirez2007,Sun2019b}.
In these models, the fractional order can be any $(0,1)$-valued function of the time variable, but in practice, it is fitted to piecewise constant time domain data, see, e.g. \cite{Adolfsson2003}. This is quite reminiscent of time-fractional 
multi-state regime-switching option pricing models. Indeed, although the CO time-fractional Black-Scholes equation can account for the nonlocal properties of the assets' prices, see e.g. \cite{Cartea2013}, it cannot describe changes in market states. Hence, regime-switching VO time-fractional models have been proposed in \cite{Chen2023,Saberi2018}
to price options. Since each market state is described by a CO time fractional model, the fractional order of the VO model is piecewise constant, see, e.g. \cite{Chen2023}[Eq. (1)].


When solving a problem with time-dependent fractional derivatives, one can distinguish two major cases:
\begin{description}
 \item [The highest order of the time derivative is a constant] This covers parabolic and hyperbolic situations; see  
 \cite{Schumer2003FractalMS,VanBockstal2022d,VanBockstal2023,Wang2019,Zheng2021b,ZHENG2022}. In this framework, the VO time derivatives can be interpreted as Volterra operators of lower order. Consequently, they can be incorporated into the right-hand side of the governing differential equation. The well-posedness of the setting can be shown using a generalised Gr\"{o}nwall lemma, see  \cite[Lemma 7.1.1]{henry} or \cite[Lemma 1]{Chen1992}.
  We refer the reader to \cite{Elliott92,slod92.1} for a discrete version of this lemma. Up to now, the most general article addressing this proof technique is \cite{slod22.2}, which also includes nonlinear functions of VO fractional derivatives.
  \item [The highest order of the time derivative is time-dependent] 
  Kochubei \cite{kochubei2011general,kochubei2019general} developed a general FC for 
\[
    D_\phi f(t):= \frac{d}{dt} \left( \phi\ast f\right)(t) - \phi(t)f(0),
\]
which is based on the Laplace transform. The integral kernel has to obey some conditions constraining $\phi$ to be completely monotone (cf. \cite{GARRAPPA2021,Luchko2020}). This is very restrictive for the choice of the kernels and, therefore, less interesting. The authors of \cite{GARRAPPA2021} present a few examples of practical relevance (exponential, Mittag-Leffler, and error function-type transitions).
Until now, there is no general theory that covers this situation. One would like to have the existence and uniqueness of a solution for $\beta(t)$, which is a piecewise smooth function. This case remains an open problem.
\end{description}

\paragraph {\bf Highlights of the paper} We consider the time-fractional diffusion problem \eqref{eq:problem} with the time-dependent highest order $\beta(t)$ of the fractional derivative in time. We assume that    $\beta(t)$ is a piecewise constant function with a finite number of steps. In Section~\ref{sec:proof_main_result}, we design a new proof technique to establish the well-posedness of the problem (stated in \Cref{thm-main}) by using Fourier analysis and the interpretation of a solution for CO FC in terms of the Mittag-Leffler function. We want to point out that $\beta(t)$ does not need to be monotone, and the convolution kernel is not positive definite. Our proof technique is limited to a finite number of steps/jumps of $\beta(t)$. Until now, we have not found a way to treat a continuously varying $\beta(t)$.

\subsection{Formulation of the problem}
Let $\Omega \subset \R^d$, $d \in \N:=\{1,2,\ldots \}$, be an open bounded domain with Lipschitz boundary $\partial \Omega$.  We consider a general second-order linear differential operator $\sL$ defined by
\bel{eq:second_order_operator}
\sL(u) = -\div \left(\mathbf{A}(x) \nabla u\right) + c(x) u, 
\ee
where 
$$
\mathbf{A}(x) = \left(a_{i,j}(x) \right)_{i,j=1,\ldots,d},\ a_{ij} \in L^\infty(\Omega,\R),\ \mathbf{A}^T = \mathbf{A},\ c \in L^\infty(\Omega,\R).
$$
Set $\sH:=L^2(\Omega)$. We assume that there exist two constants $\alpha \in (0,\infty)$ and $c_0 \in [0, \alpha \kappa^2)$, where $\kappa$ is the $\Omega$-related constant of the classical Poincar\'e inequality $\norm{\nabla u}_{\sH^d} \ge \kappa \norm{u}_{\sH}$ for all $u \in \sV:=H_0^1(\Omega)$, the closure of $C_0^\infty(\Omega)$ in the first order Sobolev space $H^1(\Omega)$, such that 
\bel{h-ell}
    \sum_{i,j=1}^d a_{ij}(x)\xi_i\xi_j \ge \alpha\abs{\boldsymbol{\xi}}^2, \quad \text{ for a.e.}\ x \in \Omega \text{ and all } \boldsymbol{\xi}=(\xi_i)_{i=1,\ldots,d} \in\mathbb{R}^d,
\ee
and
\bel{h-lb}
c(x) \ge -c_0\ \text{for a.e.}\ x \in \Omega. 
\ee
Moreover, for $T \in (0,\infty)$ fixed and for $\beta : (0,T) \to (0,1)$,  
we define the variable-order fractional integral operator ${}_0I_t^{\beta (t)} $, and the variable-order fractional Caputo operator  $  \frac {\partial ^{\beta(t)} } {\partial t ^{\beta(t)}}$ as follows  \cite{Coimbra2003,Lorenzo2002,Sun2009,Zhuang2009}
\begin{equation*}
    {}_0I_t^{\beta (t)} u(t): = \frac{1}{{\Gamma \left( {\beta (t)} \right)}}\int_0^t {\frac{{u(s)}}{{(t - s)^{1 - \beta (t)} }} \dd s} ,
\end{equation*}
\begin{equation}\label{eq:caputo_var_order}
    \partial^{\beta(t)}_t u(t):= {}_0I_t^{1 - \beta (t)} u^\prime(t) = \frac{1}{{\Gamma\left( 1 -  {\beta (t)} \right)}}\int_0^t {\frac{{u^\prime(s)}}{{(t - s)^{\beta (t)} }}\dd s}.
\end{equation}

In this contribution, we examine the existence and uniqueness issue for the solution to the following initial-boundary value problem (IBVP)
\begin{equation} \label{eq:problem}
     \left\{ \begin{array}{rlr}
    \left(\partial^{\beta(t)}_t u\right)   (x,t) +  \sL u(x,t)   &= f(x,t), & \qquad (x,t) \in \Omega \times (0,T),\\
     u(x,t) &= 0, & \qquad  (x,t) \in \partial \Omega\times (0,T), \\
     u(x,0) & = u_0(x), & \qquad  x \in \Omega, 
     \end{array} \right.
\end{equation}
in the special case where the function $\beta$ is piecewise constant. Here $u_0$ (resp., $f$) is a suitable initial condition (resp., source term) that will be made precise further.


\section{Results}
\label{sec:results}

\subsection{Definitions and notations}

In what follows, the usual norm in $\sH$ or $\sH^d$ is denoted by $\norm{\cdot}$.
We introduce the following bilinear form
$$ \ell(u,v) := \sum_{i,j=1}^d \int_{\Omega} \left( a_{i,j}(x) \partial_{x_i} u(x) {\partial_{x_j} v(x)} dx + c(x) u(x) {v(x)} \right) dx,\ u, v \in \sV. $$
Evidently, $\ell$ is continuous in $\sV \times \sV$, as we have
\begin{eqnarray*}
| \ell(u,v) | & \le & \max_{1 \le i,j \le d} \| a_{i,j} \|_{L^\infty(\Omega)} \| \nabla u \| \| \nabla v \| + \| c \|_{L^\infty(\Omega)} 
\| u \| \| v \| \\
& \le & C \| u \|_{H^1(\Omega)} \| v \|_{H^1(\Omega)},\ u,v \in \sV,
\end{eqnarray*}
where $C=\max_{1 \le i,j \le d} \| a_{i,j} \|_{L^\infty(\Omega)}+\| c \|_{L^\infty(\Omega)}$. 
Further, it readily follows from \eqref{h-ell}-\eqref{h-lb} and from the Poincar\'e inequality that $\ell$ is $\sV$-coercive:
$$
\ell(u,u)  \ge \alpha \| \nabla u \|^2 -c_0 \| u \|^2\ge (\alpha \kappa^2 -c_0) \| u \|_{H^1(\Omega)}^2,\ u \in \sV.
$$
Therefore, the linear operator $\cL$ is associated with $(\ell,\sV, \sH)$ in the sense of \cite[Chap. VI, Section 3.2.5]{Dautray1988} is selfadjoint (and positive) in $\sH$, and acts on its dense domain $D(\cL)=\{u \in \sV:\ \sL u\in \sH\}$ as $\cL u=\sL u$. 
Moreover, since $\sV$ is compactly embedded in $\sH$, the resolvent of $\cL$ is compact, and consequently, the spectrum of $\cL$ is discrete. We denote by $\lambda_n$, $n \ge 1$, the eigenvalues of $\cL$, arranged in non-decreasing order and repeated with the (finite) multiplicity (see, e.g., \cite[Theorem XIII.64]{Reed1982}):
$$ 0<  \lambda_1 \le \lambda_2 \le \cdots $$
and by $\{ X_n,\ n \ge 1 \}$ an  orthonormal basis in $\sH$ of eigenfunctions such that $\cL X_n = \lambda_n X_n$. \\
We recall that 
$$ D(\cL) = \left\{ v \in \sH : \ \sum_{n=1}^\infty \lambda_n^2 \abs{\langle v , X_n \rangle}^2 < \infty \right\}, $$
is a Banach space with respect to the norm
$$\norm{v}_{D(\cL)} := \left(\sum_{n=1}^\infty \lambda_n^2 \abs{\langle v , X_n \rangle}^2\right)^{1 \slash 2}.$$

In the next definition, we formulate how the solution to \eqref{eq:problem} should be understood. We set $I:=(0,T)$ and we recall that the space $W^{1,1}(I,\sH)$ consists of functions $u \in L^1(I,\sH)$ satisfying $\partial_t u \in L^1(I,\sH)$. 

\begin{defn}
\label{def-sol}
Let $u_0 \in D(\cL)$ and let $f \in L^1(I,\sH)$. Then, a solution to the IBVP \eqref{eq:problem} is any function
$$u \in L^1(I,D(\cL)) \cap W^{1,1}(I,\sH) $$
satisfying the two following conditions simultaneously:
\begin{enumerate}[i)]
\item $\partial_t^\beta u(x,t) + \cL u(x,t) = f(x,t)$ for a.e. $(x,t) \in \Omega \times I$,
\item $\lim_{t \downarrow 0} \norm{u(\cdot,t)-u_0}=0$.
\end{enumerate}
\end{defn}

\subsection{Main result}

Let $M \in \N$. Given $\beta_j \in (0,1)$, $j=0,1,\ldots,M-1$, and $t_0:=0<t_1<t_2<\ldots<t_M:=T$, we consider $\beta : \overline{I} \to (0,1)$ such that
$$ \beta(t):= \beta_j,\ t \in [t_j,t_{j+1}),\ j=0,1,\ldots,M-1. $$
Moreover, we put $I_j:=(t_j,t_{j+1})$, $j=0,\ldots,M-1$. 
Then, the main result of this paper can be stated as follows.

\begin{thm}
\label{thm-main}
Let $u_0 \in D(\cL)$, let $f \in W^{1,1} \left( \cup_{j=0}^{M-1} I_j,\sH \right)$, and assume that for all $j=0,\ldots,M-1$,
\bel{aa2} 
\exists \ep_j \in (0,1-\beta_j) 
\text{ such that } 
\norm{(\cdot-t_j)^{\beta_j+\epsilon_j}f^\prime}_{L^\infty(I_j,\sH)} < \infty.
\ee
Then, the IBVP \eqref{eq:problem} admits a unique solution 
$$u \in \cC^0(\overline{I},D(\cL)) \cap W^{1,1} (I,\sH) $$ 
in the sense of Definition \ref{def-sol}. Moreover, there exists a positive constant $C$, depending only on $\Omega$, $T$ and 
$\{ (t_j,\beta_j, \epsilon_j),\ j=0,\ldots,M-1 \}$, such that 
\begin{align} 
& \norm{u}_{\cC^0(\overline{I},D(\cL))} + \norm{u}_{W^{1,1}(I,\sH)} \nonumber \\
\le & C \left( \norm{u_0}_{D(\cL)} + \sum_{j=0}^{M-1} \left( 
\norm{f}_{W^{1,1}(I_j,\sH)} + \norm{(\cdot -t_j)^{\beta_j+\epsilon_j}f^\prime}_{L^\infty(I_j,\sH)} \right) \right). \label{e0} 
\end{align}
\end{thm}


\subsection{Outline and comments}

Making use of an auxiliary result stated in \Cref{pr1}, of which the proof is postponed to \Cref{sec-pr1}, we will prove the main result of this article, \Cref{thm-main}, in \Cref{sec:proof_main_result}. 

Theorem~\ref{thm-main} claims existence of a unique solution $u$ in the sense of Definition \ref{def-sol} to the IBVP \eqref{eq:problem} provided that the initial state $u_0 \in D(\cL)$ and the source term $f \in W^{1,1} \left( \cup_{j=0}^{M-1} I_j,\sH \right)$ satisfies the 
condition \eqref{aa2}. According to it, the first-order time derivative 
$f^\prime(\cdot,t)$ should not blow up faster as $t \downarrow t_j$, $j=0,\ldots,M-1$, than the power function $(t-t_j)^{-\kappa_j}$ for some $\kappa_j \in (0,1)$.
Such a condition is a key element to the proof of Theorem~\ref{thm-main}. The reason is as follows.

For all $j=0,\ldots,M-1$, $u_{| I_j}$ is characterised as a solution on the interval $I_j$ to a time-fractional PDE of constant order $\beta_j$, with source term $f_j$. Namely, we have $f_0=f$, while $f_j-f$ is expressed in terms of the first-order time derivative of the functions $u_{| I_k}$, $k=0,\ldots,j-1$, when $j=1,\ldots,M-1$, see \eqref{bb3} and \eqref{bb6}.
In the peculiar case where $j=1$ in \eqref{bb3}, we have
$$f_1^\prime(\cdot,t)=  f^\prime(\cdot,t) + \frac{\beta_1}{\Gamma(1-\beta_1)} \int_{0}^{t_1} (t-s)^{-1-\beta_1} u^\prime(\cdot,s) \dd s,\ t \in I_1, $$
so we need a good control on $u_{| I_0}^\prime$ in order to guarantee that $f_1^\prime(\cdot,t) \in L^1(I_1,\sH)$.
This is achieved by enforcing the condition \eqref{aa2} with $j=0$ on $f$, see \eqref{bb11}.

Finally, we point out that the statement of Theorem~\ref{thm-main} is no longer valid when $M$ goes to infinity. This can be understood from
the identity \eqref{def:RHS_estimates} and the estimate \eqref{e2} below (where the notation $v_j$ stands for $u_{| I_j}$), as the constant in 
the inequality \eqref{e0} blows up when $M$ becomes infinitely large.

\section{Proof of Theorem~\ref{thm-main}}
\label{sec:proof_main_result}

The proof of Theorem~\ref{thm-main} relies on the Fourier method, results of CO FC and an auxiliary result (\Cref{pr1}). We split the proof into three parts, which can be found in the next subsections.

\subsection{The Fourier method}

Using the Fourier method, we build the solution to the IBVP \eqref{eq:problem}. Namely, assuming that $u$ is a solution to \eqref{eq:problem} in the sense of Definition \eqref{def-sol}, we write
$$u(\cdot,t)=\sum_{n=1}^\infty u_n(t) X_n,\ t \in I, $$
where $u_n(t):= \langle u(\cdot,t), X_n \rangle_{\sH}$. Evidently, for each $n \in \N$, $u_n$ is the solution to the following fractional differential system (FDS)
\bel{bb1} 
\left\{ \begin{array}{ll} \pd_t^\beta u_n + \lambda_n u_n = f_{0,n}, & t \in I, \\ u_n(0)=u_{0,n}, & \end{array} \right. 
\ee
with source term
$f_{0,n}(t) :=\left\langle f(\cdot,t),X_n\right\rangle_{\sH}$ and initial state $u_{0,n}:=\left\langle u_0,X_n\right\rangle_{\sH}$. Therefore, the function $v_{j,n}:={u_n}_{| \overline{I_j}}$, $j=0,1,\ldots,M-1$, $n \in \N$, solves 
\bel{bb2}  
\left\{ \begin{array}{ll} \displaystyle  \int_{t_j}^{t} \frac{(t-s)^{-\beta_j}}{\Gamma(1-\beta_j)} v_{j,n}^\prime(s) \dd s + \lambda_n v_{j,n}(t) = f_{j,n}(t), & t \in I_j, \\
v_{j,n}(t_j) = \upsilon_{j,n}, & \end{array} \right. 
\ee
where
\bel{bb3}
f_{j,n}(t):=f_{0,n}(t) - \sum_{k=0}^{j-1} \int_{t_k}^{t_{k+1}} \frac{(t-s)^{-\beta_j}}{\Gamma(1-\beta_j)} v_{k,n}^\prime(s) \dd s,\quad t \in I_j,
\ee
and
\bel{bb4}
{\upsilon_{0,n} := u_{0,n},} \quad \upsilon_{j,n}:=v_{j-1,n}(t_j), \quad {j=1,\ldots,M-1}.
\ee
Here and in the remaining part of this text, any finite sum over $k=0$ to $j-1$ is taken equal to zero when $j=0$.

As a consequence, for all $j=0,1,\ldots,M-1$ and all $n \in \N$, we get by substituting $t-t_j$ for $t$ in \cite[Theorem~5.15]{Kilbas2006}, that the solution to \eqref{bb2} reads
\bel{bb5}
v_{j,n}(t) := \upsilon_{j,n} E_{\beta_j,1}(-\lambda_n (t-t_j)^{\beta_j}) + \int_{t_j}^{t} (t-s)^{-1+\beta_j} E_{\beta_j,\beta_j}(-\lambda_n(t-s)^{\beta_j}) 
f_{j,n}(s) \dd s,\ t \in \overline{I_j},
\ee
where $E_{\alpha,\beta}$ denotes the two-parameter Mittag-Leffler function defined by $E_{\alpha,\beta}(z) = \sum_{k=0}^\infty \frac{z^k}{\Gamma(\alpha k + \beta)}$, see \cite[Eq. (1.56)]{Podlubny1998}. By \cite[Theorem~1.6]{Podlubny1998}, we have for $\alpha \in (0,2)$ and $\beta>0$ that 
\bel{p1} 
\exists C_E>0,\ \forall z \in [0,\infty),\ \abs{E_{\alpha,\beta}(-z)} \le \frac{C_E}{1+z} \le C_E.
\ee

\subsection{Auxiliary result}
\label{sec-tech}
With reference to \eqref{bb3} and \eqref{bb5}, we set
\bel{bb6} 
v_j(\cdot,t):=\sum_{n=1}^\infty v_{j,n}(t)X_n,\quad f_j(\cdot,t):=\sum_{n=1}^\infty f_{j,n}(t) X_n,\quad t \in \overline{I_j}
, \quad \text{ for } j=0,1,\ldots,M-1,
\ee
and we note from this that 
\[
{u}_{| \overline{I_j}} = \sum_{n=1}^\infty {u_n}_{| \overline{I_j}} X_n = \sum_{n=1}^\infty v_{j,n} X_n = v_j. 
\]
The proof of \Cref{thm-main} essentially relies on the following technical result.  

\begin{prop} 
\label{pr1}
Suppose that $u_0$ and $f$ satisfy the assumptions of  \Cref{thm-main}. Then, putting
\begin{equation}\label{def:RHS_estimates}
    \mathcal{F}_j := \norm{u_0}_{D(\sL)} + \sum_{k=0}^j \left( \norm{f}_{W^{1,1}(I_k,\sH)} + \norm{(\cdot-t_k)^{\beta_k+\epsilon_k}f^\prime}_{L^\infty(I_k,\sH)} \right), \quad j=0,\ldots, M-1, 
\end{equation}
we have for all $j=0,\ldots,M-1$:
\begin{enumerate}[(i)]
\item $f_j \in W^{1,1}(I_j,\sH)$ and the estimates
\bel{e1}
\norm{f_j}_{W^{1,1}(I_j,\sH)} \le C_{j} \mathcal{F}_j,
\ee
 \begin{equation} 
\label{b0a}
\norm{f_j^\prime(\cdot,t)}  \le C_j  \mathcal{F}_j (t-t_j)^{-\beta_j-\epsilon_j}, \quad t\in I_j.
\end{equation}
\item $v_j \in \cC^0(\overline{I_j},D(\cL)) \cap W^{1,1}(I_j,\sH)$ and the estimates
\begin{equation} \label{e2}
\norm{v_j}_{\cC^0(\overline{I_j},D(\cL))} + \norm{v_j}_{W^{1,1}(I_j,\sH)} 
\le C_{j} \mathcal{F}_j,
\end{equation}
\begin{equation} \label{b0b}
\norm{v_j^\prime(\cdot,t)} \le C_{j}  \mathcal{F}_j (t-t_j)^{-1+\beta_j}, \quad t\in I_j.
\end{equation}
\end{enumerate}
Here and in the remaining part of this text, $C_j$ denotes a generic positive constant depending only on  
$\{ (I_k,\beta_k,\epsilon_k),\ k=0,\ldots,j \}$, which may change from line to line.
\end{prop}

The proof of \Cref{pr1} being quite tedious, we postpone it to \Cref{sec-pr1}.

We notice from \Cref{pr1} that $\upsilon_j \in D(\cL)$ for all $j=0,\ldots,M-1$ (this follows from the identity $\upsilon_0=u_0$ and the assumption $u_0 \in D(\cL)$ when $j=0$, and from \Cref{pr1}(ii), namely from $v_j \in \cC^0(\overline{I_j},D(\cL))$, when $j=1,\ldots,M-1$, since $v_{j-1}(\cdot,t_j)=\upsilon_j$ in this case). Moreover, we have $v_j(\cdot,t_j)=\upsilon_j$ for all $j=0,\ldots,M-1$, by virtue of \eqref{bb5}-\eqref{bb6}, and consequently
\bel{eq-pr1}
v_{j-1}(t_j)=v_j(t_j),\ j=1,\ldots,M-1.
\ee
Now, since $v_j \in \cC^0(\overline{I_j},D(\cL))$ for all $j=0,\ldots,M-1$, by \Cref{pr1}(ii), \eqref{eq-pr1} and the identity $v_j = u_{| \overline{I_j}}$ then yield that 
\bel{eq-c0}
u \in \cC^0(\overline{I},D(\cL)).
\ee

\subsection{End of the proof of Theorem~\ref{thm-main}}

Let us first prove that $u^\prime \in L^1(I,\sH)$. To this end, taking into account that $u \in L^1(I,\sH)$, according to \eqref{eq-c0} and the embedding $\cC^0(\overline{I},D(\cL)) \subset L^1(I,\sH)$, we find for every $\varphi \in \cC_0^\infty(I)$ and a.e. $x \in \Omega$, that
\beas
\langle u^\prime(x,\cdot) , \varphi \rangle_{\cC_0^\infty(I)^\prime \times \cC_0^\infty(I)} & = & 
- \langle u(x,\cdot) , \varphi^\prime \rangle_{{\cC_0^\infty(I)}^\prime \times \cC_0^\infty(I)} \nonumber \\ 
& = & -\sum_{j=0}^{M-1} \int_{t_j}^{t_{j+1}} v_j(x,t) \varphi^\prime(t) \dd t \nonumber \\
& = & \sum_{j=0}^{M-1} \int_{t_j}^{t_{j+1}} v_j^\prime(x,t) \varphi(t) \dd t - \sum_{j=0}^{M-1} \left( v_j(x,t_{j+1}) \varphi(t_{j+1}) - v_j(x,t_j) \varphi(t_j) \right). 
\eeas
Here we used that $v_j \in W^{1,1}(I_j,\sH)$ for all $j=0,\ldots,M-1$, as stated in \Cref{pr1}(ii). Further, bearing in mind that $v_j \in \cC^0(\overline{I_j},\sH)$ and that $\varphi(0)=\varphi(T)=0$, we see that
$$ \sum_{j=0}^{M-1} \left( v_j(x,t_{j+1}) \varphi(t_{j+1}) - v_j(x,t_j) \varphi(t_j) \right) = 0,\ x \in \Omega. $$
Thus, we have
$$ \langle u^\prime(x,\cdot) , \varphi \rangle_{{\cC_0^\infty(I)}^\prime, \cC_0^\infty(I)} = \sum_{j=0}^{M-1} \int_{t_j}^{t_{j+1}} v_j^\prime(x,t) \varphi(t) \dd t,\ x \in \Omega, $$
for all $\varphi \in \cC_0^\infty(I)$, from where we get that
\bel{eq-ep1} 
u^\prime(x,t) = \sum_{j=0}^{M-1} \chi_{I_j}(t) v_j^\prime(x,t),\ x \in \Omega,\ t \in I,
\ee
where $\chi_I$ denotes the characteristic function of $I$. As a consequence we have $u^\prime \in L^1(I,\sH)$ from \Cref{pr1}(ii), and hence $u \in \cC^0(\overline{I},D(\cL)) \cap W^{1,1}(I,\sH)$ by \eqref{eq-c0}.
Next, for all $j=0,\ldots,M-1$ we have
$$
\int_{t_j}^t \frac{(t-s)^{-\beta_j}}{\Gamma(1-\beta_j)} v_j^\prime(x,s) \dd s+ \cL v_j(x,t) = f_j(x,t),\ x \in \Omega,\ t \in \overline{I_j},
$$
from the first equation of \eqref{bb1} and \eqref{bb6}, whence
$$
\int_{0}^t \frac{(t-s)^{-\beta_j}}{\Gamma(1-\beta_j)} u^\prime(x,s) \dd s+ \cL u(x,t) = f(x,t),\ x \in \Omega,\ t \in \overline{I_j},
$$
by \eqref{bb3} and \eqref{eq-ep1}. Moreover, we have $u(\cdot,0)=u_0$ from \eqref{bb5}-\eqref{bb6} and the identities $u_{| \overline{I_0}} = v_0$ and $\upsilon_0=u_0$. This proves that $u$ is a solution to the IBVP \eqref{eq:problem} in the sense of Definition \ref{def-sol}. There exists at most one solution to \eqref{eq:problem} since the solution to \eqref{bb2}  on $I_j$, $j=0,\ldots,M-1$, is unique according to \cite[Theorem~4.3]{Kilbas2006}.

Finally, since $u(\cdot,t)=\sum_{j=0}^{M-1} \chi_{\overline{I_j}}(t) v_j(\cdot,t)$ for all $t \in \overline{I}$, we have
\bel{eq-ep2} 
\norm{u}_{\cC^0(\overline{I},D(\cL)} = \max \{ \norm{v_j}_{\cC^0(\overline{I_j},D(\cL)},\ j=0,\ldots,M-1 \},
\ee
and
\bel{eq-ep3}
\norm{u}_{W^{1,1}(I,\sH)} = \sum_{j=0}^{M-1} \norm{v_j}_{W^{1,1}(I_j,\sH)},
\ee
according to \eqref{eq-ep1}, so we obtain \eqref{e0} by combining \eqref{e2} with \eqref{eq-ep2}-\eqref{eq-ep3}.

\section{Proof of Proposition~\ref{pr1}}
\label{sec-pr1}

The proof is by induction on $j$. Given $j \in \{1,\ldots,M-1 \}$, the induction hypothesis (IH)$_j$ is that the claim of \Cref{pr1} (and in particular the estimates \eqref{e1} -- \eqref{b0b}) holds when $k=0,1,\ldots,j-1$ is substituted for $j$ in its statement. 

As the method of the proof is basically the same for $j=0$ and for $j \in \{1,\ldots,M-1 \}$, 
we start with the inductive step, only pointing out in the initial step the specifics of the proof for $j=0$ when needed. 

\subsection{Inductive step}
Let $j \in \{1,\ldots,M-1 \}$ be fixed. Assuming (IH)$_j$, we aim to prove the statement of \Cref{pr1}.

\subsubsection{Proof of (i)} 
We split the proof into two main steps.

\paragraph{\it Step 1}
We start by showing that $f_j \in L^1(I_j,\sH)$ satisfies the estimate
\begin{equation}
\norm{f_j}_{L^1(I_j,\sH)} 
\le C_{j}  \mathcal{F}_{j}.\label{ee1}
\end{equation}
To this purpose we refer to \eqref{bb3} and \eqref{bb6}, and obtain that
\bel{bb6b}
f_j(\cdot,t)=\sum_{n=1}^\infty \left(f_{0,n}(t) - \sum_{k=0}^{j-1} \int_{t_k}^{t_{k+1}} \frac{(t-s)^{-\beta_j}}{\Gamma(1-\beta_j)} v_{k,n}^\prime(s) \dd s\right)X_n,\ t \in I_j.
\ee
Next, for all $N \in \N$ we have
\begin{align}
 &\norm{\sum_{n=1}^N \left(f_{0,n}(t) - \sum_{k=0}^{j-1} \int_{t_k}^{t_{k+1}} \frac{(t-s)^{-\beta_j}}{\Gamma(1-\beta_j)} v_{k,n}^\prime(s) \dd s\right)X_n} \nonumber \\
 & \le \norm{\sum_{n=1}^N f_{0,n}(t)X_n} +  \frac{1}{\Gamma(1-\beta_j)}\sum_{k=0}^{j-1} \int_{t_k}^{t_{k+1}} (t-s)^{-\beta_j} \norm{\sum_{n=1}^N v_{k,n}^\prime(s)X_n}\dd s \nonumber \\
& \le \norm{f(\cdot,t)}+ \frac{1}{\Gamma(1-\beta_j)}\sum_{k=0}^{j-1} \int_{t_k}^{t_{k+1}} (t-s)^{-\beta_j} \norm{v_k^\prime(s)}\dd s, 
\nonumber
\end{align}
and $(t-s)^{-\beta_j} \le (t-t_{k+1})^{-\beta_j}$ whenever $s \le t_{k+1} \le t_j < t$ and $k=0,\ldots,j-1$. This and \eqref{bb6b} yield that 
$$
\norm{f_j(\cdot,t)} \le \norm{f(\cdot,t)}+ \frac{1}{\Gamma(1-\beta_j)} \sum_{k=0}^{j-1} (t-t_{k+1})^{-\beta_j} \norm{v_k^\prime}_{L^1(I_k,\sH)},\ t \in I_j.
$$
Thus, we get, by integrating with respect to $t$ over $I_j$ and using that
$(t_{j+1}-t_{k+1})^{1-\beta_j} - (t_{j}-t_{k+1})^{1-\beta_j} 
\le t_{j+1}^{1-\beta_j}$ for all $k=0,\ldots,j-1$, that
$$
 \norm{f_j}_{L^1(I_j,\sH)} 
 \le \norm{f}_{L^1(I_j,\sH)}+ \frac{t_{j+1}^{1-\beta_j}}{\Gamma(2-\beta_j)}\sum_{k=0}^{j-1}  \norm{v_k^\prime}_{L^1(I_k,\sH)}.
$$
Now, applying (IH)$_j$, and more precisely \eqref{e2} where $k=0,\ldots,j-1$ is substituted for $j$, we find that
\begin{equation*}
\norm{f_j}_{L^1(I_j,\sH)} \le  \norm{f}_{L^1(I_j,\sH)} +  \frac{t_{j+1}^{1-\beta_j}}{\Gamma(2-\beta_j)} \sum_{k=0}^{j-1} C_k \mathcal{F}_k,
\end{equation*}
and \eqref{ee1} follows from this and the estimates $\mathcal{F}_k \le \mathcal{F}_{j-1}$ for $k=0,\ldots,j-1$. 

\paragraph{\it Step 2} 
Having established \eqref{ee1}, we turn now to the proof of \eqref{b0a}. This can be achieved by differentiating \eqref{bb6b} with respect to $t \in I_j$: we obtain that
$$ f_j^\prime(\cdot,t)=  f^\prime(\cdot,t) + \frac{\beta_j}{\Gamma(1-\beta_j)} \sum_{k=0}^{j-1} \int_{t_k}^{t_{k+1}} (t-s)^{-1-\beta_j} v_k^\prime(\cdot,s) \dd s, $$
and hence that
\bel{bb9}
\norm{f_j^\prime(\cdot,t)}  \le \norm{f^\prime(\cdot,t)} + \frac{\beta_j}{\Gamma(1-\beta_j)} \sum_{k=0}^{j-1} \int_{t_k}^{t_{k+1}} (t-s)^{-1-\beta_j} \norm{v_k^\prime(\cdot,s)} \dd s.
\ee 
The main difficulty lies in handling the term $k=j-1$ in the above summation. We distinguish in two cases. In the first one, when $k=0,\ldots,j-2$, we get upon using that $(t-s)^{-1-\beta_j} \le (t_{j}-t_{j-1})^{-\beta_j-1}$ for all $t \in I_j$ and $s \in I_k$,  that
\begin{align}
\int_{t_k}^{t_{k+1}} (t-s)^{-1-\beta_j} \norm{v_k^\prime(\cdot,s)} \dd s 
& = \int_{t_k}^{t_{k+1}} (t-s)^{-1-\beta_j} (s-t_k)^{-1 + \beta_k} \norm{(s-t_k)^{1 - \beta_k} v_k^\prime(\cdot,s)} \dd s \nonumber
\\
& \le  (t_j-t_{j-1})^{-1-\beta_j} \left( \int_{t_k}^{t_{k+1}} (s-t_k)^{-1+\beta_k} \dd s \right) \norm{(\cdot-t_{k})^{1-\beta_k}v_k^\prime}_{L^\infty(I_k,\sH)} \nonumber \\
& \le  \beta_k^{-1} (t_j-t_{j-1})^{-1-\beta_j} (t_{k+1}-t_k)^{\beta_k} \norm{(\cdot-t_{k})^{1-\beta_k}v_k^\prime}_{L^\infty(I_k,\sH)}. \label{bb9b}
\end{align}
In the second case, when $k=j-1$, we have $(t-s)^{-\beta_j-\epsilon_j} \le (t-t_j)^{-\beta_j-\epsilon_j}$ and 
 $(t-s)^{-1+\epsilon_j} \le (t_j-s)^{-1+\epsilon_j}$ 
for all $t \in I_j$ and $s \in I_{j-1}$, whence
\begin{align*}
&\int_{t_{j-1}}^{t_j} (t-s)^{-1-\beta_j} \norm{v_{j-1}^\prime(\cdot,s)} \dd s\\
& = \int_{t_{j-1}}^{t_j} (t-s)^{-\beta_j-\epsilon_j} (t-s)^{-1+\epsilon_j}  (s-t_{j-1})^{-1+\beta_{j-1}} (s-t_{j-1})^{1-\beta_{j-1}}  \norm{v_{j-1}^\prime(s)} \dd s \\
& \le (t-t_j)^{-\beta_j-\epsilon_j}  \left( \int_{t_{j-1}}^{t_j}  (t_j-s)^{-1+\epsilon_j} (s-t_{j-1})^{-1+\beta_{j-1}}\dd s \right) \norm{(\cdot-t_{j-1})^{1-\beta_{j-1}}v_{j-1}^\prime}_{L^\infty(I_{j-1},\sH)}.
\end{align*}
Further, by performing the change of variable $r=\frac{s-t_{j-1}}{t_j-t_{j-1}}$ in the integral of the above line, we see that
\begin{align*}
\int_{t_{j-1}}^{t_j}  (t_j-s)^{-1+\epsilon_j} (s-t_{j-1})^{-1+\beta_{j-1}}\dd s
& = (t_j-t_{j-1})^{-1+\beta_{j-1}+\epsilon_j} \left( \int_0^1 (1-r)^{-1+\epsilon_j} r^{-1+\beta_{j-1}}\dd r \right) \\
& =  (t_j-t_{j-1})^{-1+\beta_{j-1}+\epsilon_j} \cB(\beta_{j-1},\epsilon_j),
\end{align*}
where $\cB$ denotes the beta function
$$ \cB(r_{1},r_{2}):=\int _{0}^{1} s^{r_{1}-1}(1-s)^{r_{2}-1} \dd s,\quad r_1 \in (0,\infty),\ r_2 \in (0,\infty). $$
As a consequence, we have for all $t \in I_j$,
\begin{align*}
& \int_{t_{j-1}}^{t_j} (t-s)^{-1-\beta_j} \norm{v_{j-1}^\prime(\cdot,s)} \dd s\\
& \le  (t_j-t_{j-1})^{-1+\beta_{j-1}+\epsilon_j} \cB(\beta_{j-1},\epsilon_j) \norm{(\cdot-t_{j-1})^{1-\beta_{j-1}}v_{j-1}^\prime}_{L^\infty(I_{j-1},\sH)}   (t-t_j)^{-\beta_j-\epsilon_j}.\end{align*}
Putting this together with \eqref{bb9}-\eqref{bb9b}, and applying \eqref{b0b} where $k=0,\ldots,j-1$ is substituted for $j$, in accordance with (IH)$_j$, we find for a.e. $t \in  I_j$ that the difference
\[(t-t_j)^{\beta_j+\epsilon_j}\norm{f_j^\prime(\cdot,t)} -
(t-t_j)^{\beta_j+\epsilon_j}  \norm{f^\prime(\cdot,t)}\]
is upper bounded by the constant 
$$ 
\frac{\beta_j}{\Gamma(1-\beta_j)} \left( \frac{(t_{j+1}-t_j)^{\beta_j+\epsilon_j}}{(t_j-t_{j-1})^{1+\beta_j}} \sum_{k=0}^{j-2} \frac{C_k(t_{k+1}-t_k)^{\beta_k}}{\beta_k} \mathcal{F}_k 
+ \frac{C_{j-1}\cB(\beta_{j-1},\epsilon_j)}{(t_j-t_{j-1})^{1-\beta_{j-1}-\epsilon_j}}   \mathcal{F}_{j-1} \right) \le C_{j} \mathcal{F}_{j-1}.
$$
Now, \eqref{b0a} follows readily from this and \eqref{aa2}. Moreover, using \eqref{ee1} and recalling that $\beta_j+\epsilon_j \in (0,1)$, we get \eqref{e1} by integrating \eqref{b0a} over $I_j$.

We turn now to showing (ii), but prior to that, we notice for further use from \eqref{e1} and from the continuous embedding $W^{1,1}(I_j,\sH) \subset \cC^0(\overline{I_j},\sH)$, that 
\bel{bb7b}
\norm{f_j(\cdot,t)} \le C_{j} \mathcal{F}_j,\quad t \in \overline{I_j}.
\ee

\subsubsection{Proof of (ii)} 
The proof of the second claim, (ii),  of \Cref{pr1} being quite lengthy, we break it into two parts.

\paragraph{\it Step 1} The first step is to show that $v_j \in \cC^0(\overline{I_j},D(\cL))$ satisfies 
$\norm{v_j}_{\cC^0(\overline{I_j},D(\cL))} \le C_{j} \mathcal{F}_j$. To do that, we recall from (i) that the function $t \mapsto f_{j,n}(t)=\langle f_j(\cdot,t) , X_n \rangle \in W^{1,1}(I_j)$, for all $n \in \N$. Then, bearing in mind that
\[
    \frac{\dd}{\dd t}   E_{\alpha,1}(-\lambda_n t^\alpha) = -\lambda_n t^{\alpha-1}E_{\alpha,\alpha}\left(-\lambda_n t^\alpha\right),\ t >0.
\] 
and integrating by parts, we obtain for all $t \in \overline{I_j}$ that
\begin{multline*}
   \int_{t_j}^t (t-s)^{-1+\beta_j} E_{\beta_j,\beta_j}(-\lambda_n(t-s)^{\beta_j})  
f_{j,n}(s) \dd s\\
 =  \frac{1}{\lambda_n} \left( E_{\beta_j,1}(0) f_{j,n}(t) - E_{\beta_j,1}(-\lambda_n (t-t_j)^{\beta_j}) f_{j,n}(t_j) + \int_{t_j}^t E_{\beta_j,1}(-\lambda_n(t-s)^{\beta_j}) f_{j,n}^\prime(s) \dd s \right). 
\end{multline*}
As a consequence, we have for all $N \in \N$ that
\begin{align}
& \norm{\sum_{n=1}^N \left( \int_{t_j}^t (t-s)^{-1+\beta_j} E_{\beta_j,\beta_j}(-\lambda_n(t-s)^{\beta_j})  
f_{j,n}(s) \dd s \right)  X_n}_{D(\cL)} \nonumber \\
& \le  \norm{\sum_{n=1}^N \left(  E_{\beta_j,1}(0) f_{j,n}(t) - E_{\beta_j,1}(-\lambda_n (t-t_j)^{\beta_j}) f_{j,n}(t_j) \right) X_n}  \nonumber \\
& \quad  + \norm{\sum_{n=1}^N \left( \int_{t_j}^t E_{\beta_j,1}(-\lambda_n(t-s)^{\beta_j}) f_{j,n}^\prime(s) \dd s \right) X_n},\quad t \in \overline{I_j}.  \label{bb8b}
\end{align}
The first term on the right-hand side of \eqref{bb8b} is simply handled by \eqref{p1}, as follows:
\begin{align}
&  \norm{\sum_{n=1}^N \left(  E_{\beta_j,1}(0) f_{j,n}(t) - E_{\beta_j,1}(-\lambda_n (t-t_j)^{\beta_j}) f_{j,n}(t_j) \right) X_n} \nonumber \\
& \le  \norm{\sum_{n=1}^N E_{\beta_j,1}(0) f_{j,n}(t) X_n} +\norm{\sum_{n=1}^N E_{\beta_j,1}(-\lambda_n (t-t_j)^{\beta_j}) f_{j,n}(t_j)  X_n}\nonumber \\
& \le  C_E \left( \left( \sum_{n=1}^N \abs{f_{j,n}(t)}^2 \right)^{1 \slash 2} + \left( \sum_{n=1}^N \abs{f_{j,n}(t_j)}^2 \right)^{1 \slash 2} \right) \nonumber \\
& \le  C_E \left( \norm{f_{j}(\cdot,t)} + \norm{f_{j}(\cdot,t_j)} \right).  \label{bb8c}
\end{align}
As for the second term, we get through standard computations that
\begin{align*}
\norm{\sum_{n=1}^N \left( \int_{t_j}^t E_{\beta_j,1}(-\lambda_n(t-s)^{\beta_j}) f_{j,n}^\prime(s) \dd s \right) X_n} 
& \le  \int_{t_j}^t \norm{\sum_{n=1}^N E_{\beta_j,1}(-\lambda_n(t-s)^{\beta_j}) f_{j,n}^\prime(s) X_n} \dd s \nonumber\\
& \le  C_E \int_{t_j}^t \left( \sum_{n=1}^N \abs{f_{j,n}^\prime(s)}^2 \right)^{1 \slash 2} \dd s \nonumber \\
& \le  C_E \int_{t_j}^t \norm{f_j^\prime(s)} \dd s.
 \end{align*}
This and \eqref{bb8b}-\eqref{bb8c} then yield that for all $t \in \overline{I_j}$,
\begin{multline} \label{bb8e}
    \norm{\sum_{n=1}^\infty \left( \int_{t_j}^t (t-s)^{-1+\beta_j} E_{\beta_j,\beta_j}(-\lambda_n(t-s)^{\beta_j}) 
f_{j,n}(s)  \dd s \right)X_n}_{D(\cL)} \\
 \le  C_E \left(  2 \norm{f_j}_{\cC^0(\overline{I_j},\sH)} + \norm{f_j^\prime}_{L^1(I_j,\sH)} \right).
\end{multline}
Now, taking into account that
\begin{align}\label{bb8}
\norm{\sum_{n=1}^\infty \upsilon_{j,n} E_{\beta_j,1}(-\lambda_n (t-t_j)^{\beta_j})X_n}_{D(\cL)}
&= \left(\sum_{n=1}^\infty \lambda_n^2 \abs{\upsilon_{j,n}}^2 E_{\beta_j,1}(-\lambda_n (t-t_j)^{\beta_j})^2\right)^{\frac{1}{2}} \nonumber \\
& \le  C_E \left(\sum_{n=1}^\infty \lambda_n^2 \abs{\upsilon_{j,n} }^2\right)^{\frac{1}{2}} \nonumber \\
& \le  C_E \norm{\upsilon_j}_{D(\cL)},\quad t \in \overline{I_j}, 
\end{align}
and recalling \eqref{bb5}-\eqref{bb6}, we deduce from \eqref{bb8e} that
\bel{bb8f}
\norm{v_j(\cdot,t)}_{D(\cL)} \le C_E \left( \norm{\upsilon_j}_{D(\cL)}+ 2 \norm{f_j}_{\cC^0(\overline{I_j},\sH)}+ \norm{f_j^\prime}_{L^1(I_j,\sH)} \right),\quad t \in \overline{I_j}.
\ee
Next, bearing in mind that $\upsilon_j = v_{j-1}(\cdot,t_j)$ since $j  \ge 1$ here, we infer from (IH)$_j$ (namely, upon substituting $j-1$ for $j$ in \eqref{e2}) that
\bel{bb8f4} 
\norm{\upsilon_j}_{D(\cL)} \le C_{j-1} \mathcal{F}_{j-1}.
\ee
By inserting this, \eqref{e1} and \eqref{bb7b} into \eqref{bb8f}, we obtain that
\bel{bb8f2}
\norm{v_j(\cdot,t)}_{D(\cL)} \le C_{j} \mathcal{F}_j,\quad t \in \overline{I_j}.
\ee
Further, we see from \eqref{bb5} that for all fixed $N \in \N$, the function $t \mapsto \sum_{n=1}^N v_{j,n}(t) X_n \in \cC^0(\overline{I_j},D(\cL))$. 
Moreover, with reference to \eqref{bb5}-\eqref{bb6}, it is clear from \eqref{bb8e}-\eqref{bb8} that the series 
$\sum_{n=1}^\infty v_{j,n}(t) X_n$ converges to $v_j(\cdot,t)$ in $D(\cL)$, uniformly in $t \in \overline{I_j}$. Therefore, we have
\bel{bb8g}
v_j \in \cC^0(\overline{I_j},D(\cL))
\ee
and the estimate
\bel{bb8f3}
\norm{v_j}_{\cC^0(\overline{I_j},D(\cL))} \le C_{j} \mathcal{F}_j,
\ee
according to \eqref{bb8f2}.

\paragraph{\it Step 2} Having established \eqref{bb8g}-\eqref{bb8f3}, we turn now to proving the estimate \eqref{b0b}.
For this purpose, we differentiate \eqref{bb5} with respect to $t \in I_j$, and obtain that
\bea
v_{j,n}^\prime(t) & = & (f_{j,n}(t_j)-\lambda_n \upsilon_{j,n}) (t-t_j)^{-1+\beta_j} E_{\beta_j,\beta_j}(-\lambda_n (t-t_j)^{\beta_j}) \nonumber \\
& & 
+ \int_{t_j}^{t} (t-s)^{-1+\beta_j} E_{\beta_j,\beta_j}(-\lambda_n (t-s)^{\beta_j}) f_{j,n}^\prime(s) \dd s,\ n \in \N.
\label{bb11}
\eea
Next, using \eqref{p1}, we get for all $N \in \N$ and all $t \in I_j$ that
\begin{align*}
& \norm{\sum_{n=1}^N \left( f_{j,n}(t_j)-\lambda_n \upsilon_{j,n}) (t-t_j)^{-1+\beta_j} E_{\beta_j,\beta_j}(-\lambda_n (t-t_j)^{\beta_j} \right) X_n} \nonumber \\
& \le   \norm{\sum_{n=1}^N  f_{j,n}(t_j) (t-t_j)^{-1+\beta_j} E_{\beta_j,\beta_j}\left(-\lambda_n (t-t_j)^{\beta_j}\right) X_n} \nonumber \\
& \quad  
+ \norm{\sum_{n=1}^N \lambda_n \upsilon_{j,n} (t-t_j)^{-1+\beta_j} E_{\beta_j,\beta_j}\left(-\lambda_n (t-t_j)^{\beta_j}\right) X_n} \nonumber \\
&\le  C_E (t-t_j)^{-1+\beta_j} \left( \left( \sum_{n=1}^N \abs{f_{j,n}(t_j)}^2 \right)^{1 \slash 2} + \left( \sum_{n=1}^N \lambda_n^2 \abs{\upsilon_{j,n}}^2 \right)^{1 \slash 2} \right) \nonumber\\
&\le  C_E (t-t_j)^{-1+\beta_j} \left(\norm{f_j(\cdot,t_j)}+\norm{\upsilon_j}_{D(\cL)} \right)
\end{align*}
and that
\begin{align*}
& \norm{\sum_{n=1}^N \left( \int_{t_j}^t (t-s)^{-1+\beta_j} E_{\beta_j,\beta_j}(-\lambda_n (t-s)^{\beta_j}) f_{j,n}^\prime(s) \dd s \right) X_n} \nonumber \\
&\le \int_{t_j}^t  (t-s)^{-1+\beta_j} \norm{\sum_{n=1}^N  E_{\beta_j,\beta_j}(-\lambda_n (t-s)^{\beta_j}) f_{j,n}^\prime(s)  X_n}\dd s \nonumber \\
&\le  C_E \int_{t_j} ^t (t-s)^{-1+\beta_j} \left( \sum_{n=1}^N \abs{f_{j,n}^{\prime}(s)}^2 \right)^{1 \slash 2} \dd s \nonumber\\
&\le  C_E \int_{t_j}^t (t-s)^{-1+\beta_j}\norm{f_j^\prime(\cdot,s)}\dd s. 
\end{align*}
Therefore, it follows from \eqref{bb5}-\eqref{bb6} and \eqref{bb11} that for all $t \in I_j$,
\bea
\norm{v_j^\prime(\cdot,t)} \le
C_E \left( (t-t_j)^{-1+\beta_j} \left(\norm{f_j(\cdot,t_j)} +\norm{\upsilon_j}_{D(\cL)} \right) + \int_{t_j}^t (t-s)^{-1+\beta_j}\norm{f_j^\prime(\cdot,s)}\dd s \right). \label{bb14}
\eea
Moreover, since $\beta_j + \epsilon_j \in (0,1)$, we apply \eqref{b0a} and get for all $t \in I_j$ that 
\begin{align*}
& \int_{t_j}^t (t-s)^{-1+\beta_j}\norm{f_j^\prime(\cdot,s)}\dd s \\
&\le   \norm{(\cdot -t_j)^{\beta_j+\epsilon_j}f_j^\prime}_{L^\infty(I_j,\sH)}  \int_{t_j}^t (t-s)^{-1+\beta_j} (s-t_j)^{-\beta_j-\epsilon_j} \dd s \\
& \le  C_j \mathcal{F}_j \cB(1-\beta_j-\epsilon_j,\beta_j) (t-t_j)^{-\epsilon_j} \\
& \le  C_j \mathcal{F}_j \cB(1-\beta_j-\epsilon_j,\beta_j) (t_{j+1}-t_j)^{1-\beta_j-\epsilon_j} (t-t_j)^{-1+\beta_j},
\end{align*}
so we obtain \eqref{b0b} by plugging the above estimate, \eqref{bb7b} and \eqref{bb8f4} into \eqref{bb14}. Finally,
\eqref{e2} follows straightforwardly from \eqref{b0b} and \eqref{bb8f3}.

\subsection{Base step}
Firstly, since $f_0=f_{| I_0}$, it is clear from the assumption $f_{| I_0} \in W^{1,1}(I_0,\sH)$ and from \eqref{aa2} that the claims of (i) in \Cref{pr1}, hold for $j=0$.

Secondly, taking $j=0$ in the derivation of \eqref{bb8f}, we get that for all $t \in I_0$,
\bel{bs1}
\norm{v_0(\cdot,t)}_{D(\cL)}  \le  C_E \left(  \norm{u_0}_{D(\cL)} + 2 \norm{f}_{\cC^0(\overline{I_0},\sH)} + \norm{f^\prime}_{L^1(I_0,\sH)} \right).
\ee
Similarly, since
$$ \int_{0}^t (t-s)^{-1+\beta_0}\norm{f^\prime(\cdot,s)}\dd s \le \cB(1-\beta_0-\epsilon_0,\beta_0) t_1^{1-\beta_0-\epsilon_0} \norm{(\cdot)^{\beta_0+\epsilon_0}f^\prime}_{L^\infty(I_j,\sH)}  t^{-1+\beta_0}, $$
for all $t \in I_0$, we find upon mimicking the proof of \eqref{bb14} that
\bel{bs2}
\norm{v_0^\prime(\cdot,t)} \le C_E \left(  \norm{u_0}_{D(\cL)} + \norm{f(\cdot,0)} + \cB(1-\beta_0-\epsilon_0,\beta_0) t_1^{-1+\beta_0} \norm{ (\cdot)^{\beta_0+\epsilon_0} f^\prime}_{L^\infty(I_0,\sH)}  \right) t^{-\epsilon_0}.
\ee
Therefore, the estimates \eqref{e2}-\eqref{b0b} for $j=0$ follow from \eqref{bs1}-\eqref{bs2} and the continuity of the embedding $W^{1,1}(I_0,\sH) \subset \cC^0(\overline{I_0},\sH)$.

This completes the proof of the proposition.

\section*{Acknowledgments}

The work of K.~Van Bockstal was supported by the Methusalem programme of Ghent University Special Research Fund (BOF) (Grant Number 01M01021).


\bibliography{refs}

\begin{thebibliography}{10}

\bibitem{Adolfsson2003}
K.~Adolfsson, M.~Enelund, and S.~Larsson.
\newblock Adaptive discretization of an integro-differential equation with a
  weakly singular convolution kernel.
\newblock {\em Comput. Methods Appl. Mech. Eng.}, 192(51-52):5285--5304, 2003.

\bibitem{almeida2019}
R.~Almeida, D.~Tavares, and D.~F.~M. Torres.
\newblock {\em The Variable-Order Fractional Calculus of Variations}.
\newblock SpringerBriefs in Applied Sciences and Technology. Springer
  International Publishing, Cham, 1st ed. 2019. edition, 2019.

\bibitem{Cartea2013}
\'{A}lvaro Cartea.
\newblock Derivatives pricing with marked point processes using tick-by-tick
  data.
\newblock {\em Quantitative Finance}, 13(1):111--123, 2013.

\bibitem{Chen1992}
C.~Chen, V.~Thom{\'e}e, and L.~Wahlbin.
\newblock Finite element approximation of a parabolic integro-differential
  equation with a weakly singular kernel.
\newblock {\em Mathematics of Computation}, 58:587--602, 1992.

\bibitem{Chen2023}
X.~Chen, X.~Gong, S.-L. Lei, and Y.~Sun.
\newblock A preconditioned iterative method for a multi-state time-fractional
  linear complementary problem in option pricing.
\newblock {\em Fractal and Fractional}, 7(4), 2023.

\bibitem{Coimbra2003}
C.~F.~M. Coimbra.
\newblock Mechanics with variable-order differential operators.
\newblock {\em Annalen der Physik}, 12(11‐12):692--703, 2003.

\bibitem{Dautray1988}
R.~Dautray and J.~L. Lions.
\newblock {\em Mathematical Analysis and Numerical Methods for Science and
  Technology}, volume Volume 2. Functional and Variational Methods.
\newblock Springer, 1988.

\bibitem{Ding2020}
W.~Ding, S.~Patnaik, S.~Sidhardh, and F.~Semperlotti.
\newblock Applications of distributed-order fractional operators: A review.
\newblock {\em Entropy}, 23(1), 2021.

\bibitem{Elliott92}
C.~Elliott and S.~Larsson.
\newblock Error estimates with smooth and nonsmooth data for a finite element
  method for the cahn-hilliard equation.
\newblock {\em Math. Comp}, 58:603--630, 1992.

\bibitem{FF2012}
S.~Fedotov and S.~Falconer.
\newblock Subdiffusive master equation with space-dependent anomalous exponent
  and structural instability.
\newblock {\em Phys. Rev. E}, 85:031132, 2012.

\bibitem{Gao2022}
Y.~Gao, D.~Yin, and B.~Zhao.
\newblock A variable-order fractional constitutive model to characterize the
  rate-dependent mechanical behavior of soft materials.
\newblock {\em Fractal and Fractional}, 6(10), 2022.

\bibitem{GARRAPPA2021}
R.~Garrappa, A.~Giusti, and F.~Mainardi.
\newblock Variable-order fractional calculus: A change of perspective.
\newblock {\em Communications in Nonlinear Science and Numerical Simulation},
  102:105904, 2021.

\bibitem{henry}
D.~Henry.
\newblock {\em Geometric theory of semilinear parabolic equations}, volume 840
  of {\em Lecture Notes in Mathematics}.
\newblock Springer-Verlag, Berlin-Heidelberg-New York, 1981.

\bibitem{K2023}
Y.~Kian.
\newblock Equivalence of definitions of solutions for some class of fractional
  diffusion equations.
\newblock {\em Math. Nachr.}, 296:5617--5645, 2023.

\bibitem{KS2023}
Y.~Kian and {\'E}.~Soccorsi.
\newblock Equivalence of definitions of solutions for some class of fractional
  diffusion equations.
\newblock {\em Inverse Probl.}, 39:125005, 2023.

\bibitem{KSY2018}
Y.~Kian, {\'E}.~Soccorsi, and M.~Yamamoto.
\newblock On time-fractional diffusion equations with space-dependent variable
  order.
\newblock {\em Ann. H. Poincar\'e}, 819:3855--3881, 2018.

\bibitem{Kilbas2006}
A.~A. {Kilbas}, H.~M. {Srivastava}, and J.~J. {Trujillo}.
\newblock {\em {Theory and applications of fractional differential equations}},
  volume 204.
\newblock Amsterdam: Elsevier, 2006.

\bibitem{kochubei2011general}
A.~N. Kochubei.
\newblock General fractional calculus, evolution equations, and renewal
  processes.
\newblock {\em Integral Equations and Operator Theory}, 71(4):583--600, 2011.

\bibitem{kochubei2019general}
A.~N. Kochubei.
\newblock General fractional calculus.
\newblock {\em Handbook of Fractional Calculus with Applications}, 1:111--126,
  2019.

\bibitem{Lorenzo2002}
C.~F. Lorenzo and T.~T. Hartley.
\newblock Variable order and distributed order fractional operators.
\newblock {\em Nonlinear dynamics}, 29(1):57--98, 2002.

\bibitem{Luchko2020}
Y.~Luchko and M.~Yamamoto.
\newblock The general fractional derivative and related fractional differential
  equations.
\newblock {\em Mathematics}, 8(12), 2020.

\bibitem{MACHADO2011}
J.~T. Machado, V.~Kiryakova, and F.~Mainardi.
\newblock Recent history of fractional calculus.
\newblock {\em Communications in Nonlinear Science and Numerical Simulation},
  16(3):1140--1153, 2011.

\bibitem{Nohel1976}
J.~Nohel and D.~Shea.
\newblock Frequency domain methods for volterra equations.
\newblock {\em Advances in Mathematics}, 22(3):278--304, 1976.

\bibitem{Podlubny1998}
I.~Podlubn\'{y}.
\newblock {\em Fractional Differential Equations: An Introduction to Fractional
  Derivatives, Fractional Differential Equations, to Methods of Their Solution
  and Some of Their Applications}.
\newblock Mathematics in Science and Engineering. Elsevier Science, 1998.

\bibitem{Podlubny1999}
I.~Podlubn\'y.
\newblock {\em Fractional Differential Equations: An Introduction to Fractional
  Derivatives, Fractional Differential Equations, to Methods of Their Solution
  and Some of Their Applications}.
\newblock Mathematics in science and engineering. Academic Press, 1999.

\bibitem{Ramirez2007}
L.~E.~S. Ramirez and C.~F.~M. Coimbra.
\newblock A variable order constitutive relation for viscoelasticity.
\newblock {\em Annalen der Physik}, 16(7‐8):543--552, 2007.

\bibitem{Reed1982}
M.~Reed and B.~Simon.
\newblock {\em Methods of modern mathematical physics. {Vol}. 4. {Analysis} of
  operators. ({Metody} sovremennoj matematicheskoj fiziki. 4. {Analiz}
  operatorov). {Transl}. from the {English} by {A}. {A}. {Pogrebkov} and {V}.
  {N}. {Sushko}}.
\newblock 1982.

\bibitem{Saberi2018}
E.~Saberi, S.~R. Hejazi, and E.~Dastranj.
\newblock A new method for option pricing via time-fractional {PDE}.
\newblock {\em Asian-Eur. J. Math.}, 11(5):15, 2018.
\newblock Id/No 1850074.

\bibitem{Schumer2003FractalMS}
R.~Schumer, D.~Benson, M.~Meerschaert, and B.~Baeumer.
\newblock Fractal mobile/immobile solute transport.
\newblock {\em Water Resources Research}, 39:1296, 2003.

\bibitem{slod92.1}
M.~Slodi\v{c}ka.
\newblock Semigroup formulation of {R}othe's method: {A}pplication to parabolic
  problems.
\newblock {\em Commentationes Mathematicae Universitatis Carolinae},
  33(2):245--260, 1992.

\bibitem{slod22.2}
M.~Slodi\v{c}ka.
\newblock Some direct and inverse source problems in nonlinear evolutionary
  {PDE}s with {V}olterra operators.
\newblock {\em Inverse Problems}, 38(12):124001, oct 2022.

\bibitem{Sun2019b}
H.~Sun, A.~Chang, Y.~Zhang, and W.~Chen.
\newblock A review on variable-order fractional differential equations:
  mathematical foundations, physical models, numerical methods and
  applications.
\newblock {\em Fractional Calculus and Applied Analysis}, 22(1):27 -- 59, 01
  Feb. 2019.

\bibitem{SunChangZhangChen2019}
H.~Sun, A.~Chang, Y.~Zhang, and W.~Chen.
\newblock A review on variable-order fractional differential equations:
  mathematical foundations, physical models, numerical methods and
  applications.
\newblock {\em Fractional Calculus and Applied Analysis}, 22(1):27--59, 2019.

\bibitem{Sun2009}
H.~G. Sun, W.~Chen, and Y.~Q. Chen.
\newblock {Variable-order fractional differential operators in anomalous
  diffusion modeling}.
\newblock {\em {Physica A: Statistical Mechanics and its Applications}},
  388(21):4586--4592, 2009.

\bibitem{VanBockstal2021}
K.~Van~Bockstal.
\newblock Existence of a unique weak solution to a non-autonomous
  time-fractional diffusion equation with space-dependent variable order.
\newblock {\em Adv. Difference Equ.}, 2021:43, 2021.
\newblock Id/No 314.

\bibitem{VanBockstal2023}
K.~Van~Bockstal, M.~A. Zaky, and A.~Hendy.
\newblock On the {R}othe-{G}alerkin spectral discretization for a class of
  variable fractional-order nonlinear wave equations.
\newblock {\em Fractional Calculus and Applied Analysis}, 2023.

\bibitem{VanBockstal2022d}
K.~{Van Bockstal}, M.~A. Zaky, and A.~S. Hendy.
\newblock On the existence and uniqueness of solutions to a nonlinear variable
  order time-fractional reaction–diffusion equation with delay.
\newblock {\em Communications in Nonlinear Science and Numerical Simulation},
  115:106755, 2022.

\bibitem{Wang2019}
H.~Wang and X.~Zheng.
\newblock Wellposedness and regularity of the variable-order time-fractional
  diffusion equations.
\newblock {\em J. Math. Anal. Appl.}, 475(2):1778--1802, 2019.

\bibitem{Zheng2021b}
X.~Zheng and H.~Wang.
\newblock A time-fractional diffusion equation with space-time dependent
  hidden-memory variable order: analysis and approximation.
\newblock {\em BIT}, 61(4):1453--1481, 2021.

\bibitem{ZHENG2022}
X.~Zheng and H.~Wang.
\newblock Analysis and discretization of a variable-order fractional wave
  equation.
\newblock {\em Communications in Nonlinear Science and Numerical Simulation},
  104:106047, 2022.

\bibitem{Zhuang2009}
P.~Zhuang, F.~Liu, V.~Anh, and I.~Turner.
\newblock Numerical methods for the variable-order fractional
  advection-diffusion equation with a nonlinear source term.
\newblock {\em SIAM Journal on Numerical Analysis}, 47(3):1760--1781, 2009.

\end{thebibliography}
\bibliographystyle{abbrv}

\end{document}